\documentclass[twocolumn,showpacs,preprintnumbers,amsmath,amssymb]{revtex4}

\usepackage[breaklinks,colorlinks=true,linkcolor=blue,citecolor=red]{hyperref}

\usepackage{graphicx}% Include figure files
\usepackage{dcolumn}% Align table columns on decimal point
\usepackage{bm}% bold math
\usepackage{color}
 \usepackage{ulem}
 
\newcommand{\blue}{ \textcolor{blue}}

\usepackage{vmargin}
\setmarginsrb{1.cm}{2.5cm}{1.cm}{2cm}{0cm}{0.5cm}{0.cm}{0.5cm}
\def\EE{{\mathbb{E}}}
\def\RR{{\mathbb{R}}}

\def\bx{{\bm x}}

\def\eps{\varepsilon}

\begin{document}

\title{Stochastic Dynamics of Incoherent Branched Flow}
\author{Josselin Garnier$^{1}$,  Antonio Picozzi$^{2}$, and Theo Torres$^{2}$}
\affiliation{$^{1}$ CMAP, CNRS, Ecole polytechnique, Institut Polytechnique de Paris, 91120 Palaiseau, France}
\affiliation{$^{2}$ Laboratoire Interdisciplinaire Carnot de Bourgogne, CNRS, Universit\'e Bourgogne Europe, Dijon, France}

\begin{abstract}
Waves propagating through weakly disordered smooth linear media undergo a universal phenomenon called branched flow. Branched flow has been observed and studied experimentally in various systems by considering coherent waves. Recent experiments have reported the observation of optical branched flow by using an incoherent light source, thus revealing the key role of coherent phase-sensitive effects in the development of incoherent branched flow. By considering the paraxial wave equation as a generic representative model, we elaborate a stochastic theory of both coherent and incoherent branched flow. We derive closed-form equations that determine the evolution of the intensity correlation function, as well as the value and the propagation distance of the maximum of the scintillation index, which characterize the dynamical formation of incoherent branched flow. We report accurate numerical simulations that are found in quantitative agreement with the theory without free parameters. 
Our theory highlights the important impact of coherence and interference on branched flow, thereby providing a framework for exploring branched flow in nonlinear media, in relation with the formation of freak waves in oceans.
\end{abstract}

\pacs{
42.25.Dd,
%Wave propagation in random media
05.40.-a.
%Fluctuation phenomena, random processes, noise, and Brownian motion
}
\maketitle

{\it Introduction.} 
Waves passing through a weakly disordered smooth medium with a correlation radius larger than the wavelength form long, narrow filaments called branches. Instead of random speckle patterns, the disordered potential focuses waves into branches that split, creating a tree-like structure, known as branched flow (BF).
Originally observed in electrons \cite{Topinka01,Aidala07,Jura07,Maryenko12,Liu13} and microwave cavities \cite{Hohmann10,Barkhofen13}, BFs have then been anticipated to occur with vastly different wavelength scales \cite{heller21}. 
%In particular, BFs may serve as a catalyst for the emergence of extreme nonlinear events \cite{green19,yuan22,jiang23branching,Mattheakis15,Mattheakis16}. 
%In this regard, they have been proposed as a focusing mechanism that could explain the formation of freak waves on the ocean \cite{white88,berry05,berry08,heller08,ying11,Degueldre16}. 
They may serve as a catalyst for the emergence of extreme nonlinear events \cite{green19,yuan22,jiang23branching,Mattheakis15,Mattheakis16}, and freak waves on the ocean \cite{white88,berry05,berry08,heller08,ying11,Degueldre16}.
%The occurrence of 
BFs have also been suggested to occur for sound waves \cite{Wolfson01}, ultra-relativistic electrons in graphene \cite{Mattheakis18}, flexural waves in elastic plates \cite{jose23branched}, while they can act as a conduit for energy transmission in scattering media \cite{brandstotter}. 
BFs have  been extended to random potentials in space and time \cite{stavina22}, to periodic potentials \cite{daza21,wagemakers25}, and even to active random walks \cite{mok23}.
More recently, BFs have been observed experimentally with optical waves propagating in soap films \cite{patsyk20}, and the control of light BF through weakly disordered media has become an important challenge \cite{brandstotter,chang24,rotter17,cao}.

The formation of BFs can be explained using geometrical optics, where local maxima of the random refraction index act as lenses, creating caustics and high wave intensities, as originally described in \cite{kulkarny82,white85}. Numerical simulations show that the scintillation index (i.e., the relative variance of the intensity fluctuations) can exceed one in such cases. Recent studies have used geometrical optics or diffraction integrals in the framework of catastrophe optics \cite{berry80,nye99} to derive the scaling behavior of BF dynamics \cite{kaplan02,metzger10} and extreme waves \cite{metzger14,pradas18uniformity}.
Actually, except for some particular theoretical studies \cite{Metzger13,Berry20}, BFs have been essentially treated in the framework of ray caustics, then disregarding coherence or interference effects \cite{heller21}. 
Along this way, experiments have been carried out essentially with coherent waves, such as coherent electron waves \cite{Topinka01}, coherent microwave \cite{Hohmann10,Barkhofen13}, or with coherent laser light \cite{patsyk20}. 
On the other hand, in recent experiments, optical BFs have been studied by using incoherent light sources \cite{patsyk22}, revealing intriguing properties about the role of coherence in the formation and the evolution of BFs, such as coherent interference between the different wave fronts and the sensitivity of BFs to the coherence of the waves.

Our aim in this Letter is to elaborate a stochastic formulation of BFs by considering an initial random wave function propagating in a random potential.
Using the paraxial wave (Schr\"odinger) equation as a representative model, we show that interference effects deeply modify the statistical properties of BFs. Employing multiscale and stochastic calculus, we derive closed-form equations that give the evolution of the intensity correlation function.
In particular we describe the evolution of the scintillation index that characterizes the dynamical formation of incoherent BFs. 
We determine that the scintillation index is a function of two dimensionless parameters that we identify and that involve the statistics of the medium and of the initial field. 
The theory of the stochastic dynamics of BFs is validated by accurate numerical simulations, which are found in quantitative agreement with the theory, without using any adjustable parameter.

{\it Model.} We consider the two-dimensional paraxial wave equation \cite{andrews2005laser,tatarskii}:
\begin{equation}
\label{eq:parax0}
i \partial_z \psi_z = - {\alpha}  \partial_x^2 \psi_z + 
V(z,x) \psi_z ,\quad z>0, \, x\in \RR,
\end{equation}
starting from $\psi_{z=0}(x)=\psi_o(x)$, where $\psi_o$ is a coherent or partially coherent field and $V$ is a smooth and slowly varying potential, which we assume to be a random process. We will denote by $\EE[\cdot]$ the expectation with respect to the distribution of this random process.

We present our work in optics as a concrete example, but the paraxial wave Eq.(\ref{eq:parax0}) is widespread in physics, making the processes discussed herein broadly applicable to various systems.
In optics, the parameter $\alpha$ and the potential $V$ are related to the index of refraction $n$ as follows: 
$\alpha= {1}/({2k_o n_o})$, $V(z,x)= { k_o (n_o^2-n^2(z,x))}{/( {2}n_o})$,
where $k_o$ is the wavenumber in free space, $n_o$ is the homogeneous  background index of refraction, and $n(z,x)$ is the spatially dependent index of refraction of the medium.

We will consider two different types of initial field.\\
{\it 1.}
We will first consider the  coherent case in which the initial field is a plane wave:
$\psi_o(x)=1$.
The measured intensity is $|\psi_z(x)|^2$, the mean intensity is $\EE\left[  |\psi_z(x)|^2 \right]$, and the scintillation index (i.e., the relative variance  of the intensity)  is 
\begin{equation}
S_z (x) = \frac{\EE\big[  |\psi_z(x)|^4 \big] - \EE\left[  |\psi_z(x)|^2 \right]^2}{  \EE\left[  |\psi_z(x)|^2  \right]^2}.
\end{equation}

\noindent
{\it 2.}
We will then consider in detail the situation in which
the initial field is a coherent or partially coherent speckled field.
We will consider the two following situations: \\
{\tt (c)} $\psi_o$ is a coherent speckled field, which will be modeled as a stationary random field with Gaussian statistics and correlation radius $\rho_o$ (the width of the field correlation function). \\
{\tt (pc)} $\psi_o$ is a partially coherent speckled field, which will be modeled as a time-dependent random field with Gaussian statistics and correlation radius $\rho_o$.\\
Such fields can be generated by passing a time-harmonic plane wave through a static (case {\tt (c)}) or rotating (case {\tt (pc)}) diffuser, which features a random arrangement of scattering centers. 
This experimental setup has attracted considerable interest  due to its ability to mimic the properties of a thermal light source \cite{crosignani71,goodman20}, 
with the added advantage of controlling the spatial and temporal coherence properties from the degree of roughness of the diffuser and its rotation speed. 
It has been used to investigate speckle phenomena \cite{goodman20}, ghost imaging \cite{shapiro08,shapiro12} and incoherent BFs \cite{patsyk22}.

We will denote by $\left<\cdot\right>$ the expectation with respect to the distribution of the initial field.
In situation {\tt (c)}, the measured intensity is $|\psi_z(x)|^2$, the mean intensity is 
$\EE\left[ \left< |\psi_z(x)|^2 \right>\right]$ and the scintillation index is 
\begin{equation}
S_z^{{\rm (c)}}(x) = \frac{\EE\big[ \left< |\psi_z(x)|^4 \right>\big] - \EE\left[ \left< |\psi_z(x)|^2 \right>\right]^2}{  \EE\left[ \left< |\psi_z(x)|^2 \right>\right]^2}.
\label{def:Sc}
\end{equation}
In situation {\tt (pc)}, assuming that the response time of the photodetector is larger than the coherence time of the field, the measured intensity is $\left<|\psi_z(x)|^2 \right>$ (the averaging $\left<\cdot\right>$ is experimentally carried out by time averaging by the detector over the multiple initial conditions generated by the rotating diffuser), the mean intensity is 
$\EE\left[ \left< |\psi_z(x)|^2 \right>\right]$,
and the scintillation index is 
\begin{equation}
S_z^{{\rm (pc)}}(x) = \frac{ \EE\big[ \left< |\psi_z(x)|^2 \right>^2\big] -\EE\left[ \left< |\psi_z(x)|^2 \right>\right]^2}{ \EE\left[ \left< |\psi_z(x)|^2 \right>\right]^2}.
\label{def:Spc}
\end{equation}

{\it Coherent initial plane wave.}
In this paragraph we assume a regime in which: i) the wavelength $\lambda = 2\pi/k_o$ is much smaller than the correlation radius $\ell_c$ of the index of refraction of the medium, ii)  the variance $\sigma_n^2$ of the index of refraction is small (hence, the variance $\sigma^2 =4 \pi^2 \sigma_n^2/\lambda^2$
of the random potential satisfies $\sigma^2 \ll 1/\lambda^2$); iii) the propagation distance is large enough so that the evolution of the variance of the intensity is of order one.

This situation has been intensively studied \cite{ishimaru78,andrews2005laser}.
By a multiscale analysis closed-form equations can be derived for the field and intensity 
correlation functions  \cite{gs14,gs23}.
These equations depend on the medium statistics via the integrated medium correlation function $\gamma$ defined by
\begin{equation}
 \label{eq:defgamma}
 \gamma(x) = \int_\RR \EE[V(0,0) V(z,x) ] dz, 
 \end{equation}
which can be written in the form:
$\gamma(x) = {\sigma^2}{\ell_c} \tilde{\gamma}\big( {x}/{\ell_c}\big)$.
 As a particular example of smooth random medium, we can consider a medium with Gaussian correlation function
$\EE \big[   {V}(0,0)
  {V}(z,x) \big] = \sigma^2 \exp\big(-(x^2+z^2)/{\ell_c^2} \big)$, so that $\tilde{\gamma}(\tilde{x})=\sqrt{\pi} \exp(-\tilde{x}^2)$.
The field correlation function is  
$\EE \big[ \psi_z \big(x+\frac{y}{2} \big) \overline{\psi_z} \big(x - \frac{y}{2} \big) \big]
=\exp\big[  z (\gamma(y)-\gamma(0))\big]$ \cite{gs14}.
This shows that the mean intensity is constant in $z$ and $x$ and that the correlation radius of the field decays as $1/\sqrt{z}$ \cite{comm1}.

We introduce two relevant parameters that will play a key role: 
$X_c=   {\sigma^{2/3} \ell_c}/{\alpha^{1/3}}$ (which is dimensionless) and $z_c=  {\ell_c}/({2\sigma^{2/3} \alpha^{2/3}})$ (which is homogeneous to a length), that will be shown to correspond to the propagation distance at which the scintillation index reaches a maximum for large values of $X_c$.
They can also be expressed as $X_c = ( 2^{4/3} \pi n_o^{1/3}) \sigma_n^{2/3}  \ell_c/\lambda  $ and $z_c=(2^{-1/3} n_o^{2/3}) \ell_c / \sigma_n^{2/3}$.
 From \cite{gs14} we find that the scintillation index does not depend on $x$:
\begin{equation}
\label{eq:expressS}
S_z  = \tilde{D}_{z/z_c}(0,0)-1,
\end{equation}
where $\tilde{D}_ {\tilde{z}}(\tilde{x},\tilde{y})$ satisfies 
\begin{align}
\partial_{\tilde{z}} \tilde{D}_{\tilde{z}} =
i X_c^{-1} \partial^2_{\tilde{x}  \tilde{y}} \tilde{D}_{\tilde{z}}
+
\frac{1}{2} 
X_c^2  \tilde{\cal U}(\tilde{x},\tilde{y})  \tilde{D}_{\tilde{z}},  \label{eq:eqtildePi:plane}
\end{align}
starting from $\tilde{D}_{\tilde{z}=0}(\tilde{x}, \tilde{y} )  = 1$,
with
$
\tilde{\cal U}(\tilde{x},\tilde{y} )=   2 \tilde{\gamma}(\tilde{x}) + 2 \tilde{\gamma}(\tilde{y}) - \tilde{\gamma}(\tilde{x}+\tilde{y}) -\tilde{\gamma}(\tilde{x}-\tilde{y}) - 2\tilde{\gamma}(0)   $
[here $\tilde{x}=x/\ell_c$ and $\tilde{y}=y/\ell_c$].
This  shows that the scintillation index   is a function of $\tilde{z}=z/z_c$ 
and $X_c$ only.
Eq.(\ref{eq:eqtildePi:plane}) can be solved by the split-step Fourier method \cite{strang}. Moreover, by expanding the solution for small $\tilde{z}$, we get $S_z \simeq [ \partial_{\tilde{x}}^4\tilde\gamma(0) / {6} ] (z/z_c)^3$ 
at leading order [with $ \partial_{\tilde{x}}^4\tilde\gamma(0) =12 \sqrt{\pi}$ for a medium with Gaussian correlation].

\begin{figure}
\begin{center}
\includegraphics[width=9cm]{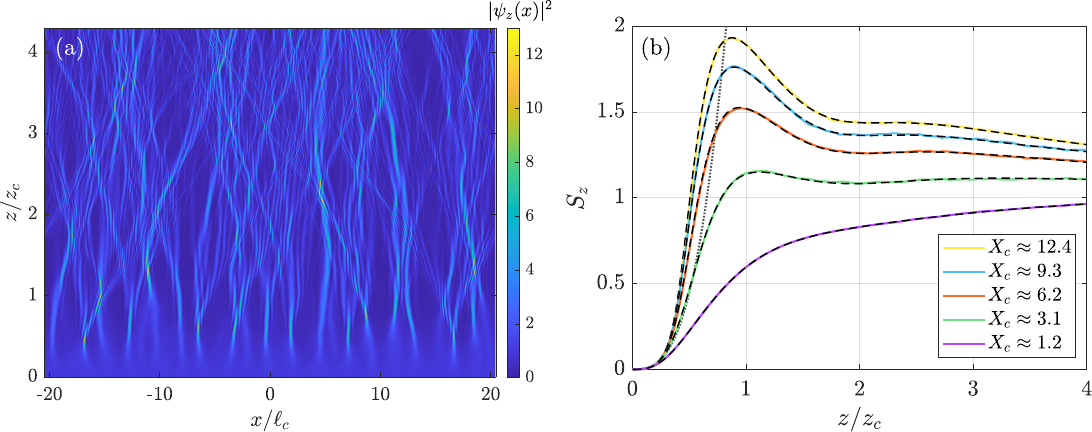} 
\end{center}
\caption{ Coherent initial plane wave  with a medium with Gaussian correlation: (a) Numerical simulation of Eq.(\ref{eq:parax0}) showing the evolution of $|\psi_z(x)|^2$ starting from  $\psi_o(x)=1$.
Parameters: $\ell_c/\lambda=100$,  
$\sigma^2 \lambda^2=10^{-4}$ ($X_c\approx 12.4$).
(b) Scintillation index $S_z$ versus $z/z_c$ for different values of $X_c$: the black dashed lines report the theory, Eq.(\ref{eq:expressS}); the dotted line is the small $z$ prediction $S_z  \simeq 2 \sqrt{\pi}  (z/z_c)^3$; the colored lines are the results of the numerical simulations, averaged over $1000$ independent realizations of the disordered potential. 
Parameters: from the bottom, $\ell_c/\lambda=10, 25, 50, 75$, with 
$\sigma^2 \lambda^2=10^{-4}$ for all curves, except for the top yellow curve ($X_c \approx 12.4$) where $\ell_c/\lambda=50$, 
$\sigma^2 \lambda^2=8 \times 10^{-4}$.
}
\label{fig:1}
\end{figure}

Here are the main results.\\
Firstly, the scintillation  index $S_z$  is close to $0$ when $z$ is small (i.e., smaller than $z_c$) and it first increases cubically with~$z$.\\
Secondly, when $X_c$ is below a threshold value $X_c^{\rm (t)}$ ($X_c^{\rm (t)}$ is between one and three for a medium with Gaussian correlation function, see Fig.~\ref{fig:1}(b)), the scintillation index is monotoneously increasing towards its limit value $1$ when $z\to +\infty$.\\
Thirdly, when $X_c$ is above the threshold value $X_c^{\rm (t)}$,
the scintillation index reaches a maximal value larger than one at finite propagation distance. The maximal value  $\max_z S_z $  depends only on $X_c$, 
 but the distance $z$ at which the maximum of the scintillation index  is reached depends on $X_c$ and $z_c$ (it was predicted to be proportional to $\ell_c / \sigma_n^{2/3}$ in previous works  \cite{kulkarny82,white85}).
 It then relaxes to its limit value $1$ when $z\to +\infty$, where the wavefield acquires Gaussian statistics for very large propagation distances \cite{metzger14,gs16,gs23,bal1,bal2}.\\
Finally, it is quite surprising to note that when $X_c$ is larger than $X_c^{\rm (t)}$, then the scintillation index may present two maxima, one global and one local [see the small bump around $z/z_c = 2.5$ in Fig.~\ref{fig:1}(b)].

We remark that the distance $\ell_c / \sigma_n^{2/3}$ is also the typical spatial scale of evolution of the number of branches that a Hamiltonian flow develops in a random potential \cite{metzger10}. As discussed in \cite[App.~A]{supplemental},  
a ray theory of branched flow
can predict the numbers and positions of local intensity maxima, however it cannot predict the values of the maxima that result from interference effects and that depend on the coherence properties of the initial field (see also \cite{Berry20}).

\begin{figure}
\begin{center}
\includegraphics[width=9cm]{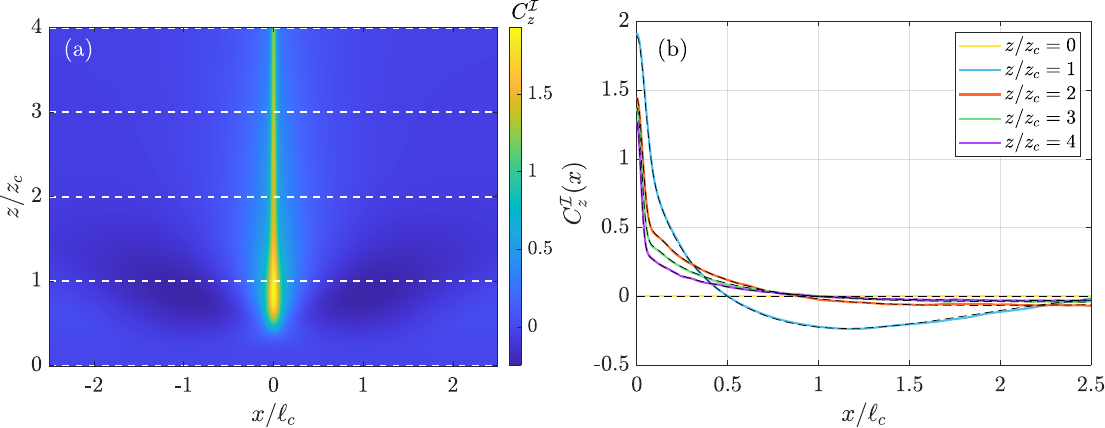} 
\end{center}
\caption{Coherent initial plane wave with a medium with Gaussian correlation: (a) Theoretical intensity correlation function $C^{{\cal I}}_z(x)$ from Eq.(\ref{eq:C_I_coh}).
(b) Comparison of $C^{{\cal I}}_z(x)$ from Eq.(\ref{eq:C_I_coh}) (black dashed lines), with the numerical simulations of Eq.(\ref{eq:parax0}) (colored lines), for different propagation lengths $z/z_c$. An average over $1000$ simulations with different realizations of the random potential $V(z,x)$ has been carried out.
Parameters: $\ell_c/\lambda=100$, $\sigma^2 \lambda^2=10^{-4}$
($X_c \approx 12.4$).
}
\label{fig:2}
\end{figure}

We have also computed the intensity correlation function
\begin{equation}
C^{\cal I}_z(x) = 
\frac{
\EE\big[ 
\big| \psi_z \big(y+ \frac{x}{2}\big) \big|^2
\big| \psi_z \big(y-\frac{x}{2}\big) \big|^2 \big]-
\EE\big[ 
\big| \psi_z \big(y\big) \big|^2\big]^2}{
\EE\big[ 
\big| \psi_z \big(y\big) \big|^2 \big]^2},
\end{equation}
which is independent of $y$ and is given by
\begin{equation}
C^{\cal I}_z(x) =   \tilde{D}_{z/z_c}(x/\ell_c,0)-1,
\label{eq:C_I_coh}
\end{equation}
where $\tilde{D}$ is solution of Eq.(\ref{eq:eqtildePi:plane}).
The intensity correlation function is plotted   in Fig.~\ref{fig:2}. Of course one has $C^{\cal I}_{z=0}(x)=0$, $C^{\cal I}_z(0)=S_z$,
$ C^{\cal I}_z(x)\to 0$ as $x\to +\infty$,
and   $\int C^{\cal I}_z(x) dx=0$ (this can be interpreted as an energy conservation relation).

We have tested the validity of the theoretical predictions by direct numerical simulations of the paraxial wave Eq.(\ref{eq:parax0}) (see \cite[App.~H]{supplemental}). The results for the evolution of the scintillation index in Fig.~\ref{fig:1} and the intensity correlation function in Fig.~\ref{fig:2} show excellent quantitative agreements, even though the separation of scales is not strong in the simulations.

 {\it Incoherent initial wave: Scaling regime.}
From now on we address the situation in which the initial field $\psi_o(x)$ is a speckled field.
In addition to the assumptions (i) to (iii) considered above for the initial plane-wave case, we assume that 
%the initial field 
$\psi_o(x)$ has Gaussian statistics, with 
%a 
correlation radius $\rho_o$ 
%that is 
larger than the wavelength $\lambda$ and smaller than the correlation radius $\ell_c$ of the index of refraction, $\ell_c/\lambda \gg \rho_o/\lambda \gg 1$.

We carry out a multiscale analysis in which a dimensionless scale parameter $\eps$ encapsulates the four assumptions listed above.
Accordingly, 
we denote by $\eps \sim \lambda/\ell_c$, the order of magnitude of the ratio of the wavelength over the correlation radius of the index of refraction. 
We assume that the typical amplitude of the fluctuations of the index of refraction is $\eps^c$, with $c>0$.
If we consider that the reference length is the correlation radius of the index of refraction, 
we can write $\alpha^\eps=1/(2k_on_o) = \eps \alpha$ and $V^\eps =k_o(n_o^2-n^2)/( {2} n_o )= \eps^{c-1} {V}$, and the scaled paraxial wave equation has the form
\begin{equation}
\label{ew:scaledschro}
i \partial_z \psi^\eps_z = - \eps {\alpha} \partial_x^2 \psi^\eps_z +
\eps^{c-1} {V}(z,x) \psi^\eps_z  ,\quad z>0, \, x\in \RR ,
\end{equation}
starting from $\psi^\eps_{z=0}(x)=\psi_o^\eps(x)$. 
The initial field $\psi_o^\eps$ has a correlation radius of the order of $\eps^d$ (relative to the correlation radius of the index of refraction) for some $d\in (0,1)$, which means that it is larger than the wavelength (because $d<1$) and smaller than the correlation radius of the index of refraction (because $d>0$).
The correlation function of the initial field is, therefore, of the form
\begin{equation}
\left< \psi_o^\eps(x+\eps^d \frac{y}{2})\overline{\psi_o^\eps}(x-\eps^d \frac{y}{2})\right> = {\cal C}_o\big( y \big)  .
\end{equation}
For the numerical simulations
we consider the Gaussian model \cite{foley,friberg} in which ${\cal C}_o(y) = \exp ( - y^2/(4\rho_o^2))$.\\
Finally, we consider the Wigner transform after a propagation distance of order $\eps^{-b}$ (relative to the correlation radius of the index of refraction):
\begin{equation}
\label{eq:defWeps}
W^\eps_z(x,k) = \int_\RR \left< \psi^\eps_{\frac{z}{\eps^b}}(x+\eps^d \frac{y}{2})\overline{\psi^\eps_{\frac{z}{\eps^b}}}  (x-\eps^d \frac{y}{2})\right> e^{-i ky }dy   .
\end{equation}
In the scaling regime $d \in (1/5,1)$, $b=1-d$, $c=3(1-d)/2$, we get from (\ref{ew:scaledschro}) (see \cite[App.~B]{supplemental}) that it  satisfies the scaled Vlasov-type equation
\begin{equation}
\partial_z W^\eps_z +  \partial_k \omega_k \partial_x W^\eps_z
-  \frac{1}{\eps^{b/2}} \partial_x {V}\big( \frac{z}{\eps^b},x\big) \partial_k  W^\eps_z = 0 ,
\label{eq:vlasoveps}
\end{equation}
with the initial condition $ W^\eps_{z=0}(x,k) = {\cal W}_o(k )=   \int_\RR {\cal C}_o(y) e^{-iky} dy $ and with $\omega_k = \alpha k^2 $ (see \cite[App.~C]{supplemental}).
Note that the scaling of the potential  in (\ref{eq:vlasoveps}) 
is appropriate for the use of limit theorems for random differential equations  \cite[Chapter 6]{book07} and we will carry out such a multiscale analysis.

Before going to the multiscale analysis, we remark that the solution of the Vlasov equation (\ref{eq:vlasoveps}) can be expressed in terms of the solutions of random ordinary differential equations.
Indeed, using the characteristic method, we have
${W}^\eps_z  \big(X^\eps_z(x,k), K^\eps_z(x,k) \big)={\cal W}_o(k)
,
$
where $(X^\eps_z(x,k),K^\eps_z(x,k))$ satisfies the ray equations
\begin{equation}
\label{eq:difXeps}
\frac{dX^\eps_z}{dz} = 2\alpha K^\eps_z, \qquad \frac{dK^\eps_z}{dz} = - \frac{1}{\eps^{b/2}} \partial_x {V}\big(\frac{z}{\eps^b},X^\eps_z) , 
\end{equation}
starting from  $X^\eps_{z=0}(x,k)=x$, $ K^\eps_{z=0}(x,k)=k$.
A key result (proved in \cite[App.~C]{supplemental}) that makes it possible to study the Wigner transform is the following one:
For any ${X} , {K} \in \RR$, 
\begin{equation}
\label{eq:expressWeps}
{W}^\eps_z({X},{K})
=
\int_{\RR^2} {\cal W}_o(k) \delta(X^\eps_z(x,k)-{X})\delta(K^\eps_z( x,k)-{K})  dx dk  .
\end{equation}
By taking an expectation (with respect to the distribution of the random medium), one can see that the mean Wigner transform involves the probability density function (pdf) of $(X^\eps_z(x,k),K^\eps_z(x,k))$. Higher-order moments of the Wigner transform involve multivariate pdf. Those pdf are computed in \cite[App.~D]{supplemental}, and they give the following results.

\begin{figure}
\begin{center}
\includegraphics[width=9cm]{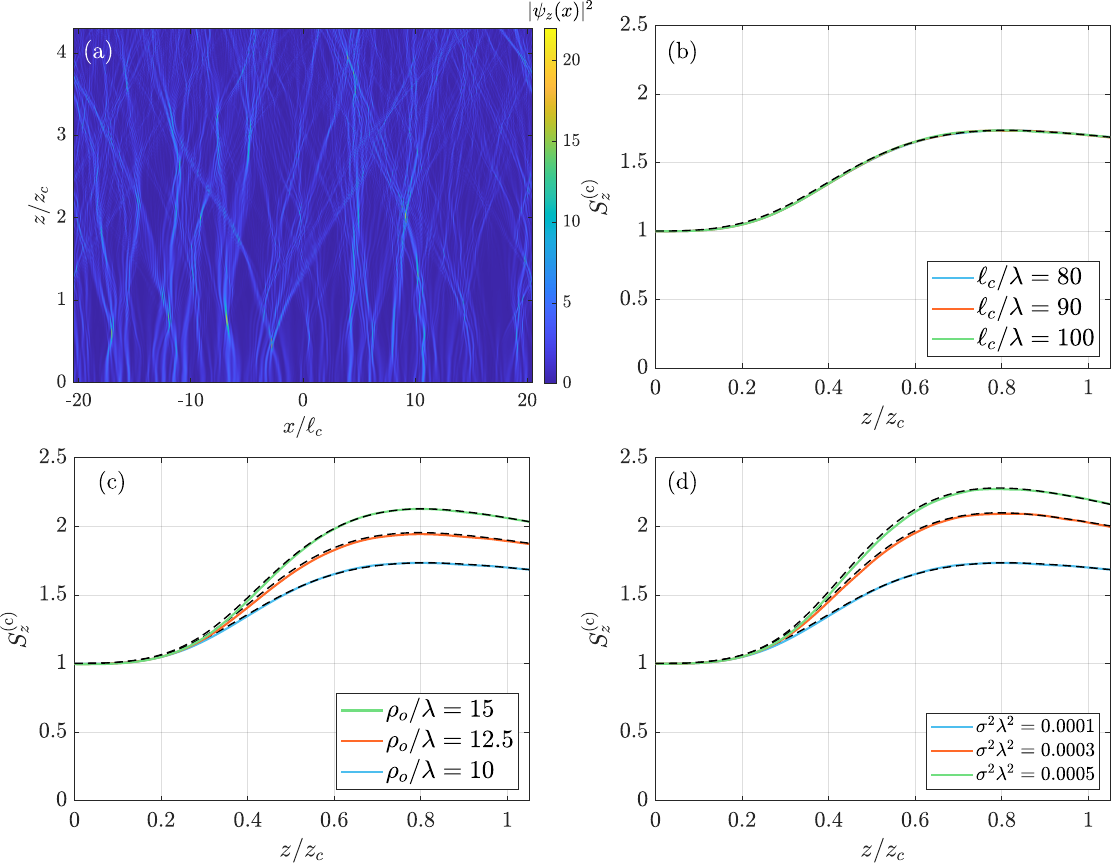} 
\end{center}
\caption{ Incoherent initial wave  with a medium with Gaussian correlation:
(a) Numerical simulation of Eq.(\ref{eq:parax0}) showing the evolution of $|\psi_z(x)|^2$ starting from a coherent speckle field [situation {\tt (c)}], with $\rho_o/\lambda=10$, $\ell_c/\lambda=100$, $\sigma^2 \lambda^2=10^{-4}$.
(b-d)
Evolution of 
%the scintillation index 
$S_z^{(c)}$ versus $z/z_c$, by varying different parameters: the black dashed lines report the theory, Eq.(\ref{eq:Szii1});  the colored lines are the results of the numerical simulations, averaged over $1000$ independent realizations of the disordered potential and of the initial random field. 
Parameters: 
(b) $\rho_o/\lambda=10$, $\sigma^2 \lambda^2=10^{-4}$; 
(c) $\ell_c/\lambda=100$, $\sigma^2 \lambda^2=10^{-4}$; 
(d) 
$\rho_o/\lambda=10$, $\ell_c/\lambda=100$.
} 
\label{fig:3}
\end{figure}

\begin{figure}
\begin{center}
\includegraphics[width=9cm]{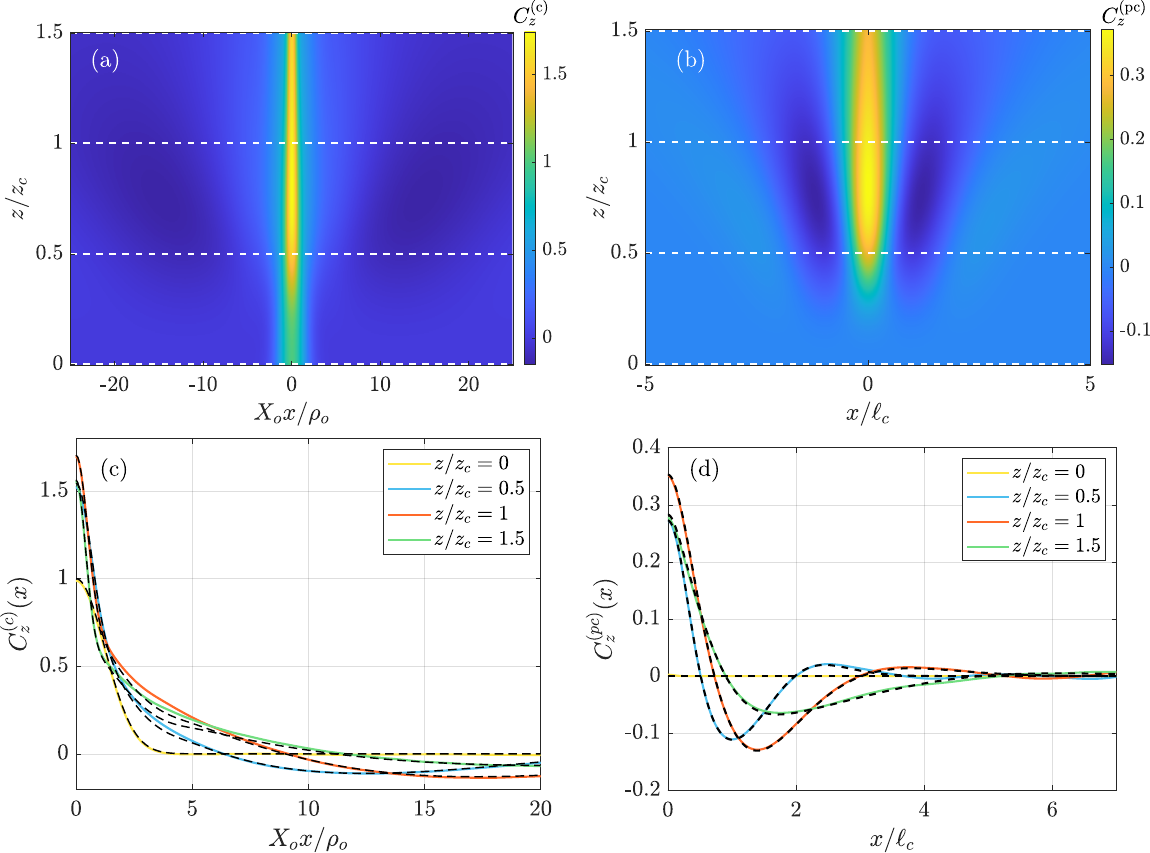} 
\end{center}
\caption{ Incoherent initial wave with a medium with Gaussian correlation:
Theoretical intensity correlation function $C^{{\cal I}}_z(x)$: from Eq.(\ref{eq:C_I_c}) for a coherent speckle field (a) [situation {\tt (c)}], and from Eq.(\ref{eq:C_I_pc}) for a partially coherent speckle field (b) [situation {\tt (pc)}].
Corresponding comparison of the theoretical correlation function $C^{{\cal I}}_z(x)$ (black dashed lines), with the numerical simulations of Eq.(\ref{eq:parax0}) (colored lines), for different propagation lengths $z/z_c$. 
In situation {\tt (c)}, an average over $1000$ 
%independent 
realizations of 
$V(z,x)$ and of 
$\psi_o(x)$, has been considered. In situation {\tt (pc)} an average over  $300$ %independent 
realizations of 
$V(z,x)$, each with $400$ 
%independent 
realizations of 
$\psi_o(x)$.
Parameters: $\ell_c/\lambda=100$, $\rho_o/\lambda=10$, 
$\sigma^2 \lambda^2=10^{-4}$.
} 
\label{fig:4}
\end{figure}

{\it Mean Wigner transform.}
From (\ref{eq:expressWeps}) we get the expression of the mean Wigner transform in the regime $\eps \to 0$,
which in turn gives  the expression of the field correlation function 
$\EE \big[ \big< \psi_z \big(x+\frac{y}{2} \big) \overline{\psi_z} \big(x - \frac{y}{2} \big) \big> \big]
= {\cal C}_o(y) \exp ( -{\gamma_2  z} y^2/2 )$, where $\gamma_2=-\partial_x^2 \gamma(0)$.
This shows that  the mean intensity is constant in $z$ and $x$ and that 
the correlation radius of the beam decays as $1/\sqrt{z}$ just as in the case of an initial coherent plane wave.

{\it Scintillation index.}
We write the correlation function of the initial field 
in the dimensionless form ${\cal C}_o(y) = \tilde{\cal C}_o( y/\rho_o)$, where $\rho_o$ is the correlation radius of the initial field.
We introduce the relevant dimensionless parameter \mbox{$X_o= {\sigma^{2/3} \rho_o}/ {\alpha^{1/3}}$}.
We get that in the situation {\tt (pc)} and {\tt (c)} the scintillation index  does not depend on $x$ (see \cite[App.~E]{supplemental}):
\begin{equation}
\label{eq:Szii1}
S_z^{{\rm (pc)}}  = \tilde{\Pi}_{z/z_c}(  0 , 0) -
1  ,\qquad
S_z^{{\rm (c)}} = 2  \tilde{\Pi}_{z/z_c}(  0 , 0) 
-
1  ,
\end{equation}
where ${\tilde{\Pi}}_{\tilde{z}}(\tilde{x},\tilde{y})$ is the solution to
\begin{align}
\partial_{\tilde{z}}  {\tilde{\Pi}}_{\tilde{z}} =
 i \partial_{\tilde{x}}\partial_{\tilde{y}}  {\tilde{\Pi}}_{\tilde{z}}
- \frac{1}{2} \big(\tilde\Gamma(0)-\tilde\Gamma(\tilde{x})\big) \tilde{y}^2 {\tilde{\Pi}}_{\tilde{z}},
\label{eq:hattildePi}
\end{align}
starting from $ \tilde{\Pi}_{\tilde{z}=0}( \tilde{x}, \tilde{y} )  =
 \tilde{\pi}_o(\tilde{y} / X_o)$. 
Here $\tilde{\Gamma}(\tilde{x}) = - \partial_{\tilde{x}}^2 \tilde{\gamma}(\tilde{x})$, $\ {\tilde{\pi}}_o(\tilde{y})= |\tilde{\cal C}_o(\tilde{y}) |^2/ \tilde{\cal C}_o(0)^2$,
while  Eq.(\ref{eq:hattildePi}) can be solved by a split-step Fourier method \cite{strang}.
By expanding the solution of  (\ref{eq:hattildePi}) for small $\tilde{z}$, we get 
 $S_z^{{\rm (pc)}} \simeq 
[\partial_{\tilde{x}}^4\tilde\gamma(0) /6]
(z/z_c)^3$, and $S_z^{{\rm (c)}} \simeq 1+ 
[\partial_{\tilde{x}}^4\tilde\gamma(0) /3]
 (z/z_c)^3$ at leading order, see \cite[App.~F]{supplemental}  [with $ \partial_{\tilde{x}}^4\tilde\gamma(0) =12 \sqrt{\pi}$ for a medium with Gaussian correlation].
This shows that the early dynamics of the scintillation index in situation {\tt (pc)} does not depend on the correlation radius $\rho_o$ nor on the correlation function of the initial field, and is equivalent to the behavior valid for an initial plane wave.
The scintillation index first grows cubically and then reaches a maximum value,
which depends on 
$X_o= {\sigma^{2/3} \rho_o}/ {\alpha^{1/3}}$
only (it increases with $X_o$).
It is interesting to note that $X_o$ (and hence the maximal scintillation indices) depends on $\rho_o$ but not on $\ell_c$, while the distance $z$ at which the maximum of the scintillation index  is reached depends on $\ell_c$ through $z_c=\ell_c / (2\sigma^{2/3}\alpha^{2/3})$.
Remember that when the initial field is a plane wave, the maximum of the scintillation index depends only on $X_c=  \sigma^{2/3} \ell_c / \alpha^{1/3}$. This shows that speckled beams experience reduced intensity growth compared to plane waves: the smaller the correlation radius of the initial beam, the lower the maximal intensity reached by the beam as it propagates.

The theory (obtained in the limit $\eps \to 0$) is compared to simulations of Eq.(\ref{eq:parax0}).
The intensity evolution in Fig.~\ref{fig:3}(a) exhibits distinct qualitative features 
with respect to the coherent excitation in Fig.~\ref{fig:1}(a).
As observed experimentally \cite{patsyk22}, in the coherent case, each branch is accompanied by sidelobes arising from interference effects, which tend to disappear when the initial condition is incoherent.
In the simulations we study the impact of $\rho_o$, $\sigma^2$ and $\ell_c$, on the scintillation index.
Fig.~\ref{fig:3} shows an excellent quantitative agreement,  in spite of the rather limited separation of scales 
of the parameters.
The simulations also confirm that the maximum of 
$S_z$ does not depend on 
$\ell_c$, see Fig.~3(b).

{\it Intensity correlation function.}
Our theoretical approach can be exploited to compute an explicit form for the fourth-order moments of the field. In particular, the intensity correlation function in situation {\tt (c)} is 
\begin{align}
C^{{\cal I},{\rm (c)}}_z(x) 
 &=  {\tilde{\Pi}}_{z/z_c}(x/\ell_c,0)  +
 {\tilde{\Pi}}_{z/z_c} (x/\ell_c,X_o x / \rho_o)  -1 .
\label{eq:C_I_c}
\end{align}
This expresses a two-scale behavior:
 At the small scale \mbox{$x \sim \rho_o$}, the intensity correlation function decays rapidly, while its behaviour exhibits complex variations at the large scale $x \sim \ell_c$. Similarly, the intensity correlation function in situation {\tt (pc)}
is  given by 
\begin{equation}
C^{{\cal I},{\rm (pc)}}_z(x) = 
 {\tilde{\Pi}}_{z/z_c} (x/\ell_c,0)  -1 .
\label{eq:C_I_pc}
\end{equation}
Note that it is equal to $C^{{\cal I},{\rm (c)}}_z (x)$ when $x$ is of the order of $\ell_c$ because $\ell_c \gg \rho_o$ and ${\tilde{\Pi}}_{\tilde{z}} (\tilde{x},\tilde{y}) \to 0$ as $\tilde{y}\to+\infty$.
The intensity correlation functions are plotted in Fig.~\ref{fig:4}.
The  two-scale behavior of the intensity correlation function in situation {\tt (c)} is clearly visible: the limit for large $x/\rho_o$ of the intensity correlation function is the initial value for $C^{{\cal I},{\rm (pc)}}_z$ (or $C^{{\cal I},{\rm (c)}}_z$) when $x$ is of the order of $\ell_c$.
This behavior can be seen in the numerical simulations as well in Fig.~\ref{fig:4},  and it satisfies the energy conservation relation $\int C^{{\cal I},{\rm (pc)}}_z(x) dx=0$.

{\it Perspectives.}  
We have reported a general stochastic theory of BFs by considering both coherent and incoherent initial
waves. Optics naturally provides an ideal framework for experimentally
testing and observing the theoretical predictions.
The results presented in the two-dimensional framework can be extended to the three-dimensional one (see \cite[App.~G]{supplemental}).
Our work paves the way for a systematic approach to studying coherent phase-sensitive effects in BFs, focusing on the linear propagation regime. As shown in prior studies, linear BFs can trigger extreme nonlinear events \cite{green19,yuan22,jiang23branching,Mattheakis15,Mattheakis16,dudley14,safari17,dudley19}, with intensity peaks significantly influenced by nonlinearity \cite{green19,jiang23}. Our stochastic framework offers a basis for developing a theoretical model of nonlinear branched flow.

\vglue.2cm

{\bf Acknowledgements.} Funding was provided by Agence Nationale de la Recherche (Grants No. ANR-23-CE30-0021, ANR-21-ESRE-0040) and by Agence de l'Innovation de D\'efense (AID) via Centre Interdisciplinaire d'\'Etudes pour la D\'efense et la S\'ecurit\'e (project PRODIPO).  Calculations were performed using HPC resources from DNUM-CCUB (Universit\'e de Bourgogne Europe).

%\newpage

\newpage

\appendix

\centerline{\bf Supplementary Material on the Article}

\section{Ray theory of branched flow}
One of the main tools used to describe and understand the properties of branched flow has been the ray tracing method.
Rays are constructed as the characteristic curves of the eikonal equation obtained by considering rapidly oscillating solution of the paraxial Eq.(1) (main text) \cite{synge:SM}
\begin{equation}
    \psi_z(x) = A(z,x)e^{iS(z,x)/\delta},
\end{equation}
where $\delta$ is an order counting parameter that expresses the assumption that the local phase $S(z,x)$ of the wave varies rapidly in comparison with its
amplitude $A(z,x)$.
Expanding the wave
equation in a hierarchical fashion in powers of $\delta^{-1}$ yields
a system of equations for the local phase and amplitude as
asymptotic series in $\delta$.
The leading-order equation is the eikonal equation: 
\begin{equation}
    -\partial_z S_0 = \alpha \left(\partial_x S_0\right)^2 + V(z,x),
\end{equation}
with $S_0$ the first term in the asymptotic series of $S$.
The eikonal equation can be interpreted as a Hamilton-Jacobi equation for a Hamiltonian $H(z,x,k_z,k_x)$ that determines the ray trajectories. The Hamiltonian is obtained by making the substitution $(\partial_zS_0,\partial_xS_0) \rightarrow (k_z,k_x)$:\begin{equation}
    H = k_z + \alpha k_x^2 + V(z,x).
\end{equation}
The condition $H=0$ is the local dispersion relation.

\begin{figure}[!h]
    \centering
    \includegraphics[width = 0.8\linewidth]{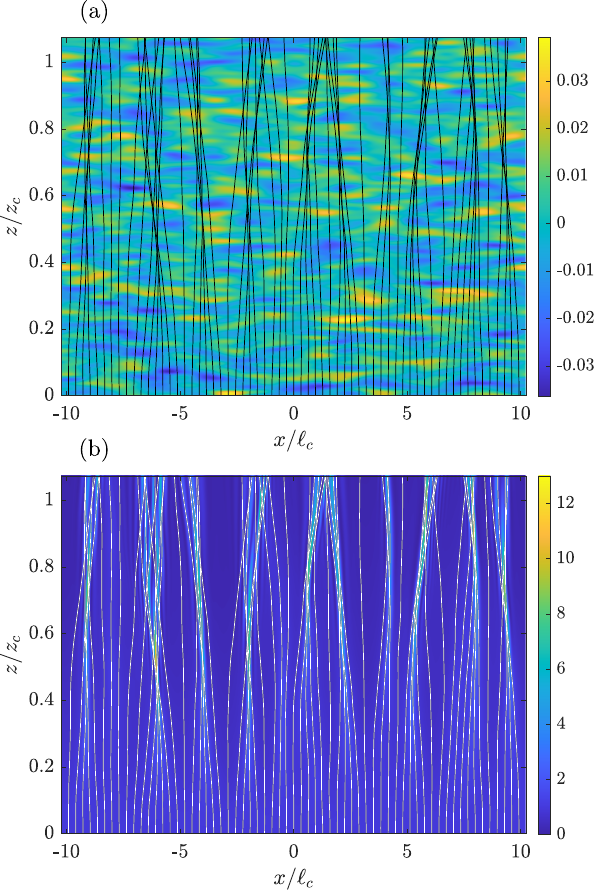}
    \caption{Coherent initial plane wave. Example of a congruence of rays for a single realization of the random potential. Panel (a) shows a realization of the random potential $V(z,x)$ together with the rays associated with an initial coherent plane wave at $z=0$. Panel (b) shows the BF obtained by solving the paraxial Eq.(1) (main text) with the potential shown in (a), superimposed with the congruence of rays. We can see the formation of caustics which are associated with increases in wave intensity. Parameters are the same as in Fig.~1(a) (main text).}
    \label{fig:supplemental1}
\end{figure}

\begin{figure}[!h]
    \centering
    \includegraphics[width = 0.8\linewidth]{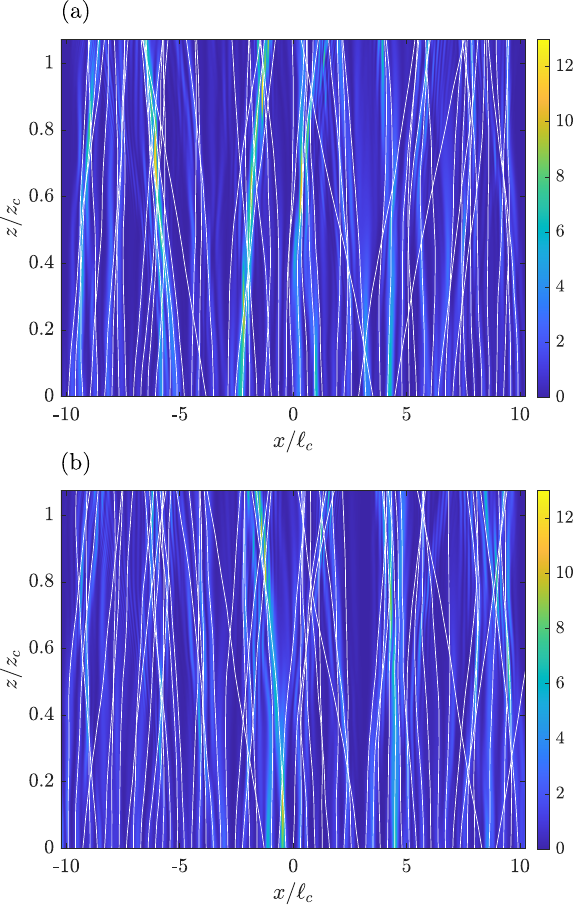}
     \caption{Incoherent initial wave. (a) Evolution of the BF obtained by solving the paraxial Eq.(1) (main text), superimposed with the congruence of rays. 
(b) Same as in (a), except that a different realization of the initial incoherent speckle field is considered.
The realization of the random potential $V(z,x)$ in (a) and (b) is the same as in Fig.~\ref{fig:supplemental1}.
As in the case of the initial plane wave, we can see the formation of caustics despite the fact that the rays are not all launched parallel to each other due to the initial phase variations. Parameters are the same as in Fig.~3(a) (main text) }
    \label{fig:supplemental2}
\end{figure}

The rays are parametrized curves, $(x(s),z(s))$, that are the
solutions of Hamilton's equations:
\begin{equation}
    \dot{x} = \partial_{k_x} H , \quad \dot{z} = \partial_{k_z} H ,\quad
    \dot{k}_x = - \partial_x H , \quad \dot{k}_z = - \partial_z H  ,
\end{equation}
where the dot represents a derivative with respect to the ray parameter $s$.
Since the Hamiltonian is linear in $k_z$, we have $\dot{z} = 1$ and therefore we can choose the $z$ coordinate to parametrize the rays.
Plugging the expression for the Hamiltonian, the equations to be solved are 
\begin{equation}
    \dot{x} = 2 \alpha k_x, \quad 
    \dot{k}_x = -  \partial_x V , \quad \dot{k}_z = - \partial_z V,
\end{equation}
together with the initial conditions $(x(z=0),k_x(z=0),k_z(z=0))$.
From the definition of the momentum, we have $k_x(z=0) = \partial_x S(z=0,x)$ and $k_z$ is then found by solving the condition $H=0$.

In the case of a coherent plane wave, $S(z=0,x)$ is constant and therefore $k_x(z=0)=0$ for all rays (all the rays are parallel to each other).
In the case of an incoherent initial speckle field, $k_x \neq 0$ and the launched rays are not parallel.
Figures~\ref{fig:supplemental1}-\ref{fig:supplemental2} report the ray dynamics and the solution of the paraxial Eq.(1) (main text) for the same realization of the random potential and for different initial conditions: coherent plane wave in Fig.~\ref{fig:supplemental1}, and incoherent speckled field in Fig.~\ref{fig:supplemental2} (panels (a) and (b) report two different realizations of the speckle field).
We can clearly see the formation of caustics which are associated with increases in wave intensities.
This illustrates the fact that the rays equations can predict the positions of the maximal intensities, but the values of the maxima result from interference effects that depend on the coherence properties of the initial field. This was discussed in \cite{Berry20:SM}  by using a simple random phase model. 
The determination of the statistics of the values of the intensity maxima require a detailed multiscale analysis as carried out in the main text.

\section{Scaled Vlasov equation}
\label{app:vla}%
We consider the Wigner transform $W^\eps_z(x,k) $  defined by \blue{Eq.(12)} (main text).
From \blue{Eq.(10)} (main text),
it  satisfies the scaled Vlasov-type equation 
\begin{align*}
& \partial_z W^\eps_z +\eps^{1-d-b}  \partial_k \omega_k \partial_x W^\eps_z  \\
& + \eps^{c-b-1} i \int_\RR \big[ {\cal V}( \frac{z}{\eps^b},x+\eps^d \frac{y}{2}) 
-{\cal V} ( \frac{z}{\eps^b},x-\eps^d \frac{y}{2}) \big] \\
&\quad \times
 \left< \psi^\eps_{ \frac{z}{\eps^b} } (x+\eps^d \frac{y}{2})\overline{\psi^\eps_{\frac{z}{\eps^b} }} (x-\eps^d \frac{y}{2})\right> \exp( -i ky ) dy = 0 ,
\end{align*}
with $\omega_k = \alpha k^2 $,
which gives after expansion of the last term of the left-hand side
\begin{align}
\nonumber
&\partial_z W^\eps_z +\eps^{1-d-b}  \partial_k \omega_k \partial_x W^\eps_z 
-  \eps^{c+d-b-1} \partial_x {\cal V}\big( \frac{z}{\eps^b},x\big) \partial_k  W^\eps_z\\
& = O(\eps^{c+3d-b-1}) ,
\label{eq:Wepsreste}
\end{align}
with the initial condition $ W^\eps_{z=0}(x,k) = {\cal W}_o( k ) $.
In the scaling regime
$d \in (1/5,1)$, $
b=1-d$, $c=3(1-d)/2$,
we have $1-d-b=0$, $c+d-b-1=-b/2$, and $c+3d-b-1 = (5d-1)/2>0$, so that we can neglect the remainder in (\ref{eq:Wepsreste}) and we get that $W^\eps_z$ satisfies
\blue{(13)} (main text)
with the initial condition $W^\eps_{z=0} (x,k) = {\cal W}_o( k ) $. 
 
 Note that the scaling regime addressed here is different from the one used to derive the paraxial white-noise (or It\^o-Schr\"odinger) model
\cite{andrews2005laser:SM,gs09a:SM,gs09b:SM,gs16:SM}. The paraxial white-noise model  is valid when $d=0$, $b=1$, $c=3/2$,  that is to say, when the wavelength is much smaller than the correlation radius of the medium, which is itself of the same order as 
the correlation radius of the initial field.

\section{Proof of Equation (15) (main text)}
\label{app:char}%
We have 
$$
{W}^\eps_z ({X},{K})
=
\int_{\RR^2} {W}^\eps_z(x',k') \delta(x'-{X})\delta(k'-{K}) dx' dk'  .
$$
We make the change of variables $(x',k') \mapsto (x,k)$ with 
$x'= X^\eps_z(x,k)$, $k'=K^\eps_z(x,k)$:
\begin{align*}
{W}^\eps_z({X},{K})
&=
\int_{\RR^2} {W}^\eps_z(X^\eps_z(x,k),K^\eps_z(x,k)) \delta(X^\eps_z(x,k)-{X}) \\
&\quad \times \delta(K^\eps_z(x,k)-{K})  
|{\rm Det} {\bf J}_z^\eps (x,k) |dx dk ,
\end{align*}
where ${\bf J}_z^\eps(x,k)$ is the Jacobian
$$
{\bf J}^\eps_z (x,k)= \begin{pmatrix}
\frac{\partial X^\eps_z}{\partial x} 
(x,k)
&
\frac{\partial X^\eps_z}{\partial k} 
(x,k)\\
\frac{\partial K^\eps_z}{\partial x} 
(x,k)
&
\frac{\partial K^\eps_z}{\partial k} 
(x,k)
\end{pmatrix}  .
$$
On the one hand we have ${W}^\eps_z(X^\eps_z(x,k),K^\eps_z(x,k)) = {\cal W}_o(k)$ and on the other hand we can compute
\begin{align*}
&\frac{d}{dz} \frac{\partial X^\eps_z}{\partial x}  = 2\alpha \frac{\partial K^\eps_z}{\partial x}  , \qquad \frac{\partial X^\eps_z}{\partial x}\mid_{z=0} (x,k) = 1, \\
&\frac{d}{dz} \frac{\partial K^\eps_z}{\partial x}  = 
-\frac{1}{\eps^{b/2}} \partial^2_x {V}(\frac{z}{\eps^b},X^\eps_z)
\frac{\partial X^\eps_z}{\partial x}   , \qquad \frac{\partial K^\eps_z}{\partial x} \mid_{z=0}(x,k) = 0,\\
&\frac{d}{dz} \frac{\partial X^\eps_z}{\partial k}  = 2\alpha \frac{\partial K^\eps_z}{\partial k}  , \qquad \frac{\partial X^\eps_z}{\partial k} \mid_{z=0}(x,k) = 0,\\
&\frac{d}{dz} \frac{\partial K^\eps_z}{\partial k}  = 
-\frac{1}{\eps^{b/2}} \partial^2_x {V}(\frac{z}{\eps^b},X^\eps_z)
\frac{\partial X^\eps_z}{\partial k}   , \qquad \frac{\partial K^\eps_z}{\partial k} \mid_{z=0}(x,k) = 1,
\end{align*}
which gives
\begin{align*}
\frac{d}{dz} {\rm Det} {\bf J}^\eps_z  &=
\Big( \frac{d}{dz} \frac{\partial X^\eps_z}{\partial x} \Big)
 \frac{\partial K^\eps_z}{\partial k}
 +
 \frac{\partial X^\eps_z}{\partial x} \Big(\frac{d}{dz}  \frac{\partial K^\eps_z}{\partial k} \Big)
\\
&\quad  -
 \Big(\frac{d}{dz}  \frac{\partial K^\eps_z}{\partial x} \Big)
  \frac{\partial X^\eps_z}{\partial k}
  -
   \frac{\partial K^\eps_z}{\partial x}
   \Big(\frac{d}{dz}  \frac{\partial X^\eps_z}{\partial k} \Big)\\
   &= 
2\alpha \frac{\partial K^\eps_z}{\partial x}  
 \frac{\partial K^\eps_z}{\partial k}
 -\frac{1}{\eps^{b/2}} \partial_x^2 {\cal V}(\frac{z}{\eps^b},X^\eps_z)
 \frac{\partial X^\eps_z}{\partial x}  \frac{\partial X^\eps_z}{\partial k}\\
 &\quad 
 +\frac{1}{\eps^{b/2}} \partial_x^2 {\cal V}(\frac{z}{\eps^b},X^\eps_z)
 \frac{\partial X^\eps_z}{\partial x}
  \frac{\partial X^\eps_z}{\partial k}
  -
2\alpha   \frac{\partial K^\eps_z}{\partial x}
  \frac{\partial K^\eps_z}{\partial k}
   \\
   &=0 ,
\end{align*}
hence ${\rm Det} {\bf J}^\eps_z  = {\rm Det} {\bf J}^\eps_{z=0} ={\rm Det}  {\bf I}=1$.
This gives the desired result \blue{Eq.(15)} (main text).

\section{Diffusion approximation theory}
\label{sec:adif}%
This section contains the technical results that are needed to characterize the statistics of the wave field, in particular the width of the envelope, the correlation radius of the field and the scintillation index.
The potential ${V}$ is a smooth, stationary, random process 
with mean zero and integrable covariance function.
Applying diffusion-approximation theory \cite[Chapter 6]{book07:SM}, we can show from \blue{(14)} (main text) that,
for any integer $n$, for any $x_1,\ldots,x_n\in \RR$, for any $k_1,\ldots,k_n \in \RR$, the $\RR^{2n}$-valued process $(X^\eps_z(x_j,k_j),K^\eps_z(x_j,k_j))_{j=1}^n$
converges in distribution as $\eps \to 0$ to the Markov diffusion process 
$({X}_z(x_j,k_j),{K}_z(x_j,k_j))_{j=1}^n$
with the infinitesimal generator 
\begin{equation}
{\cal L}^{(n)} = \sum_{j=1}^n 2\alpha {K}_j \frac{\partial}{ \partial X_j} 
+ \frac{1}{2} \sum_{j,j'=1}^n \Gamma(X_j-X_{j'}) \frac{\partial^2}{\partial K_j \partial K_{j'}}   \,,
\end{equation}
where 
\begin{equation}
\label{eq:defGamma}
\Gamma(x) = \int_{-\infty}^\infty \EE \big[ \partial_x {V}(0,0)
\partial_x {V}(z,x) \big]dz  .
\end{equation}
As a particular example of smooth random medium, we can consider a potential ${V}$ with Gaussian correlation function, variance $\sigma^2$ and correlation radius $\ell_c$.
We then have
\begin{equation}
\Gamma(x) = 2\sqrt{\pi} \sigma^2 \ell_c^{-1} \Big(1-\frac{2x^2}{\ell_c^2}\Big) \exp\Big( - \frac{x^2}{\ell_c^2}\Big).
\end{equation}

{\bf Application $n=1$.}
Let $x,k \in \RR$.
The pdf $p_z^{(1)}(X,K;x,k)$ of $({X}_z(x,k),{K}_z(x,k))$ satisfies the Fokker-Planck equation
\begin{equation}
\label{eq:eqp1z:0}
\partial_z p_z^{(1)} = ({\cal L}^{(1)})^* p_z^{(1)} , 
\end{equation}
starting from $p_{z=0}^{(1)}(X,K;x,k) = \delta(X-x) \delta(K-k)$, 
where 
$$
{\cal L}^{(1)} = 2\alpha K \partial_X   +\frac{\Gamma(0)}{2}\partial_K^2 ,
$$
and $({\cal L}^{(1)})^*$ is the adjoint of ${\cal L}^{(1)}$.
Eq.~(\ref{eq:eqp1z:0})  has the form
\begin{equation}
\label{eq:eqp1z}
\partial_z p_z^{(1)} =-2\alpha K \partial_X p_z^{(1)} +\frac{\Gamma(0)}{2}\partial_K^2 p_z^{(1)} .
\end{equation}
It is possible to solve this equation (by taking a Fourier transform in $(X,K)$) and we get the expression of the pdf of the limit process $(X_z(x,k),K_z(x,k))$:
\begin{align}
\nonumber
p_z^{(1)}(X,K;x,k) =& \frac{1}{\sqrt{2\pi\Gamma(0)z}}
\exp\Big( - \frac{(K-k)^2}{2\Gamma(0) z}\Big)
 \frac{1}{\sqrt{2\pi\Gamma(0)\frac{z^3}{3}}} \\
&\times \exp\Big( - \frac{3(X-x-\alpha(K+k)z)^2}{2\Gamma(0) z^3}\Big) .
\label{eq:pz1}
\end{align}

{\bf Application $n=2$.}
Let $x_1,x_2,k_1,k_2 \in \RR$.
The pdf $p_z^{(1)}(X_1,X_2,K_1,K_2;x_1,x_2,k_1,k_2)$ of $({X}_z(x_j,k_j),{K}_z(x_j,k_j))_{j=1}^2$
 satisfies the Fokker-Planck equation
\begin{equation}
\label{eq:eqp2z:0}
\partial_z p_z^{(2)} = ({\cal L}^{(2)})^* p_z^{(2)} , 
\end{equation}
starting from $p_{z=0}^{(2)}(X_1,X_2,K_1,K_2;x_1,x_2,k_1,k_2) = \delta(X_1-x_1) \delta(X_2-x_2) \delta(K_1-k_1)\delta(K_2-k_2)$, 
where ${\cal L}^{(2)}$ is the infinitesimal generator of $({X}_z(x_j,k_j),{K}_z(x_j,k_j))_{j=1}^2$:
\begin{align}
\nonumber
{\cal L}^{(2)} =&  2\alpha {K}_1 \frac{\partial}{ \partial X_1} +
2\alpha {K}_2 \frac{\partial}{ \partial X_2} 
+ \frac{1}{2}\Gamma(0) 
\Big( 
\frac{\partial^2}{\partial K_1^2} 
+
\frac{\partial^2}{\partial K_2^2} 
\Big)\\
&
+
\Gamma(X_1-X_2) 
\frac{\partial^2}{\partial K_1 \partial K_2}  .
\label{eq:defL2}
\end{align}
We introduce
\begin{align}
\label{eq:cv1:sm}
& R=\frac{X_1+X_2}{2}, \quad Q = X_1-X_2, \\
& U = \frac{K_1+K_2}{2},\quad V=K_1-K_2,
\label{eq:cv2:sm}
\end{align}
where $X_j=X_z(x_j,k_j)$, $K_j =K_z(x_j,k_j)$, $j=1,2$.
The infinitesimal generator of the process $(R_z,Q_z,U_z,V_z)$ is
\begin{equation}
\label{eq:defL3}
{\cal L} = 2\alpha U\partial_R +2\alpha V\partial_Q
+ \frac{1}{4} \big(\Gamma(0)+\Gamma(Q)\big) \partial_U^2 +\big(\Gamma(0)-\Gamma(Q)\big) \partial_V^2  .
\end{equation}
In particular, the process $(Q_z,V_z)$ is Markov with generator
\begin{equation}
{\cal L} =  2\alpha V \partial_Q
+ \big(\Gamma(0)-\Gamma(Q)\big) \partial_V^2  .
\label{eq:genQV:sm}
\end{equation}

\section{Expression of the scintillation index for incoherent initial conditions}
From  \blue{Eq.(15)} (main text) we get the expression of the second-order moment of the Wigner transform in the limit $\eps \to 0$:
\begin{align}
\nonumber
&
\lim_{\eps \to 0}
\EE \big[ {W}^\eps_z( {X}_1,\ {K}_1) {W}^\eps_z( {X}_2, {K}_2) \big]\\
\nonumber
&
=\int_{\RR^4} {\cal W}_o(k_1)  {\cal W}_o(k_2)  \\
&\quad \times p_z^{(2)}({X}_1,{K}_1,{X}_2,{K}_2;x_1,k_1,x_2,k_2)   dx_1 dk_1 dx_2 dk_2  ,
\label{eq:mom20:sm}
\end{align}
where $p_z^{(2)}$ is the solution of the Fokker-Planck equation (\ref{eq:eqp2z:0}).
The second-order moment of the intensity in situation {\tt (pc)} is
\begin{align}
\nonumber
&\EE \Big[ \Big< \big|\psi^\eps_{\frac{z}{\eps^b}}( {X}) \big|^2 \Big>^2  \Big]\\
&
=\frac{1}{(2\pi)^2}
 \int_{\RR^2} \EE \big[ {W}^\eps_z({X},{K}_1) {W}^\eps_z({X},{K}_2) \big]
d{K}_1d{K}_2  .
\label{eq:44:sm}
\end{align}
The second-order moment of the intensity in situation {\tt (c)}  is
\begin{align}
\nonumber
&
\EE \Big[ \Big< \big|\psi^\eps_{\frac{z}{\eps^b}}(  {X} ) \big|^4 \Big>  \Big] \\
&
=
\frac{2}{(2\pi)^2}
 \int_{\RR^2} \EE \big[ {W}^\eps_z({X},{K}_1) {W}^\eps_z({X},{K}_2) \big]
d{K}_1d{K}_2   ,
\label{eq:45:sm}
\end{align}
where we have used Isserlis' theorem \cite{isserlis18:SM} 
\begin{align*}
\left< \psi_o^\eps(x)\overline{\psi_o^\eps}(y)\psi_o^\eps(x')\overline{\psi_o^\eps}(y') \right>
=&
\left< \psi_o^\eps(x)\overline{\psi_o^\eps}(y)\right>\left<\psi_o^\eps(x')\overline{\psi_o^\eps}(y') \right> \\
&+
\left< \psi_o^\eps(x) \overline{\psi_o^\eps}(y') \right>\left<\psi_o^\eps(x')\overline{\psi_o^\eps}(y)\right>.
\end{align*}
By Eq.(\ref{eq:mom20:sm}) and the change of variables (\ref{eq:cv1:sm}-\ref{eq:cv2:sm}) 
we then get:
\begin{align}
\nonumber
&
\lim_{\eps \to 0}
\EE \Big[ \Big< \big| \psi^\eps_{\frac{z}{\eps^b}}  (  {R} )\big|^2 \Big>^2  \Big]  \\
\nonumber
&
=\frac{1}{(2\pi)^2}
 \int_{\RR^6}
 {\cal W}_o(u+\frac{v}{2})
 {\cal W}_o(u-\frac{v}{2})\\
 &\quad \times 
 p_z( {R}, 0, {U},{{V}} | r,q,u,v) 
 d{U} d{{V}} 
 dr dq du dv
 \nonumber  \\
 \nonumber
 & = \frac{1}{(2\pi)^2}
 \int_\RR 
 \Big[ \int_\RR   
 {\cal W}_o(u+\frac{v}{2})
 {\cal W}_o(u-\frac{v}{2}) d u\Big]\\
 &\quad \times
 \Big [\int_{\RR^2} p_z(0,  {{V}} | q,v) 
dq d{{V}}  \Big] dv  ,
\label{eq:fourthmoment1}
\end{align}
which does not depend on ${R}$.
By (\ref{eq:genQV:sm}),
in the last line $p_z(  {Q}, {{V}} | q,v) $ is the pdf solution of 
\begin{equation}
\partial_z p_z =
 - 2\alpha V \partial_Q p_z
+ \big(\Gamma(0)-\Gamma(Q)\big) \partial_V^2 p_z,
\end{equation}
starting from
$p_{z=0}( {Q}, {{V}} | q,v)  = \delta({Q}-q)\delta({{V}}-v)$.
Eq.(\ref{eq:fourthmoment1}) can be rewritten as
\begin{align}
&
\lim_{\eps \to 0}
\EE \Big[ \Big< \big| \psi^\eps_{\frac{z}{\eps^b}}  ( {R} )\big|^2 \Big>^2  \Big]  
=
{\Pi}_z(0,0) 
\end{align}
in terms of 
the function ${\Pi}_z$ defined by
\begin{equation}
{\Pi}_z(Q,S)= \int_{\RR^3}p_z(Q,{{V}} |q,v)   \pi_o(v) 
e^{i S V} dq dv  dV ,
\end{equation}
with $\pi_o(v) =  \frac{1}{(2\pi)^2}
\int_\RR   
 {\cal W}_o(u+\frac{v}{2})
 {\cal W}_o(u-\frac{v}{2}) d u$.
 The function ${\Pi}_z$  is the solution of
\begin{equation}
\partial_z {\Pi}_z =
 2 i \alpha   \partial_{Q}\partial_{S} {\Pi}_z
-  \big(\Gamma(0)-\Gamma(Q)\big) {S}^2 {\Pi}_z,
\label{eq:edpPiz}
\end{equation}
starting from ${\Pi}_{z=0}( {Q}, {S})  =  \int \pi_o(v) e^{iSv} dv = |{\cal C}_o({S})|^2$. This gives \blue{Eqs.(16-17)} (main text).

\section{Proof of the small $z$-expansion}
\label{app:smallZ}%
Let $\tilde{\Pi}_{\tilde{z}}$ be the solution of \blue{(17)} (main text).
We consider the functions
$$
 \tilde{M}_{j,{\tilde{z}}}(\tilde{x}) = (-i)^j 
 \partial_{\tilde{y}}^j \tilde{\Pi}_{\tilde{z}} (\tilde{x},\tilde{y}) \mid_{\tilde{y}=0} .
 $$
 They satisfy the equations 
 \begin{align*}
  \partial_{\tilde{z}} \tilde{M}_{0,{\tilde{z}}} &= - \partial_{\tilde{x}} \tilde{M}_{1,{\tilde{z}}}, \\
  \partial_{\tilde{z}} \tilde{M}_{1,{\tilde{z}}} &= - \partial_{\tilde{x}} \tilde{M}_{2,{\tilde{z}}},\\
  \partial_{\tilde{z}} \tilde{M}_{2,{\tilde{z}}} &= - \partial_{\tilde{x}} \tilde{M}_{3,{\tilde{z}}} +\big(\tilde{\Gamma}(0)-\tilde{\Gamma}(\tilde{x})\big) \tilde{M}_{0,{\tilde{z}}} ,
  \end{align*}
  starting from $ \tilde{M}_{j,{\tilde{z}}=0}(\tilde{x})=\tilde{M}_{j,o}:= (-iX_o)^j \tilde{\pi}_o^{(j)}( 0)$.
  For small ${\tilde{z}}$ and using the fact that $\tilde{M}_{j,o} $ does not depend on $\tilde{x}$,  we get successively:
\begin{align*}
\tilde{M}_{2,{\tilde{z}}} (\tilde{x})&= \tilde{M}_{2,o} +\big(\tilde{\Gamma}(0)-\tilde{\Gamma}(\tilde{x})\big) \tilde{M}_{0,o} {\tilde{z}} +o({\tilde{z}}),\\
\tilde{M}_{1,{\tilde{z}}} (\tilde{x}) &= \tilde{M}_{1,o} + \frac{1}{2} \partial_{\tilde{x}} \tilde{\Gamma}(\tilde{x}) \tilde{M}_{0,o} {\tilde{z}}^2 +o({\tilde{z}}^2),\\
\tilde{M}_{0,{\tilde{z}}} (\tilde{x}) &= \tilde{M}_{0,o} - \frac{1}{6} \partial_{\tilde{x}}^2 \tilde{\Gamma}(\tilde{x}) \tilde{M}_{0,o} {\tilde{z}}^3 +o({\tilde{z}}^3).
\end{align*}
By \blue{Eq.(16)} (main text)
this gives the desired result for the small $z$-expansions of $S^{\rm (c)}_z$ and $S^{\rm (pc)}_z$ since $\tilde{M}_{0,o}=1$.
More specifically, we get 
$S_z^{{\rm (pc)}}=  \frac{\tilde{\gamma}_4}{6} \frac{z^3}{z_c^3} + o\big(\frac{z^3}{z_c^3}\big)$, $S_z^{{\rm (c)}}= 1+  \frac{\tilde{\gamma}_4}{3}\frac{z^3}{z_c^3} + o\big(\frac{z^3}{z_c^3}\big)$, with $\tilde{\gamma}_4=-\partial_{\tilde{x}}^2 \tilde{\Gamma}(0) = \partial_{\tilde{x}}^4 \tilde{\gamma}(0)$, where $\gamma(x) = \int_\RR \EE[V(0,0) V(z,x) ] dz={\sigma^2}{\ell_c} \tilde{\gamma}\big( {x}/{\ell_c}\big)$. For a medium with Gaussian correlation function, $\tilde{\gamma}({\tilde x})=\sqrt{\pi}\exp(-{\tilde x}^2)$ and $\tilde{\gamma}_4=12 \sqrt{\pi}$.

\section{Extension to the three-dimensional case}
\label{app:3D}%
The results described in this paper can be readily extended to the three-dimensional paraxial wave equation:
\begin{equation}
\label{eq:parax3d}
i \partial_z \psi_z = - {\alpha}  \big( \partial_{x_1}^2+ \partial_{x_2}^2\big) \psi_z + 
V(z,\bx) \psi  , 
\end{equation}
for $z>0$, $\bx= (x_1,x_2)\in \RR^2$.
As an illustration, 
let us assume that that the initial field is Gaussian with Wigner transform independent of $\bx$:
\begin{equation}
\int_{\RR^2} \left< \psi_o(\bx+\frac{{\bm y}}{2})\overline{\psi_o}(\bx-\frac{{\bm y}}{2})\right> e^{ -i {\bm k} \cdot{\bm y} } d{\bm y} = \tilde{\cal W}_o(\rho_o {\bm k}),
\end{equation}
and that the medium fluctuations have Gaussian covariance function:
\begin{equation}
\label{eq:gauscov:3d}
\EE \big[   {V}(0,{\bf 0})
  {V}(z,\bx) \big] = \sigma^2 \exp\Big(-\frac{|\bx|^2+z^2}{\ell_c^2} \Big).
\end{equation}
Under such circumstances, 
in the situation {\tt (pc)} the scintillation index defined by \blue{(4)} (main text) has the form
\begin{equation}
\label{eq:Szii1:3d}
S_z^{{\rm (pc)}}  =  \tilde{\Pi}_{z/z_c}(  {\bf 0} , {\bf 0} )   
-
1  ,
\end{equation}
while in the situation {\tt (c)}  the scintillation index defined by \blue{(3)} (main text) has the form
\begin{equation}
\label{eq:Szi1:3d}
S_z^{{\rm (c)}}  =2 \tilde{\Pi}_{z/z_c}(  {\bf 0} ,{\bf 0} ) 
-
1  .
\end{equation}
The scintillation index in situations {\tt (c)} and {\tt (pc)} depends on the function $\tilde{\Pi}_{\tilde{z}}$ that is the solution of:
\begin{equation}
\partial_{\tilde{z}} \tilde{\Pi}_{\tilde{z}} =
 i  \nabla_{\tilde{\bm x}}  \cdot \nabla_{\tilde{\bm y}}\tilde{\Pi}_{\tilde{z}}
- \frac{1}{2} \sum_{j,l=1}^2 \big( \tilde{\Gamma}_{jl}({\bf 0})- \tilde{\Gamma}_{jl}(\tilde{\bm x})\big)  \tilde{y}_j \tilde{y}_l \tilde{\Pi}_{\tilde{z}}, 
\label{eq:edpPiz:3d}
\end{equation}
starting from $\tilde{\Pi}_{\tilde{z}=0}( \tilde{\bm x}, \tilde{\bm y} )  = |\tilde{\cal C}_o(\tilde{\bm y} /X_o)|^2/\tilde{\cal C}_o({\bf 0})^2$,
where $\tilde{\cal C}_o$ is the inverse Fourier transform of $\tilde{\cal W}_o$ and
\begin{align}
\tilde{\boldsymbol{\Gamma}}(\tilde{\bm x}) 
&=
2\sqrt{\pi} \bigg(
{\bf I} - 2 \begin{pmatrix}
\tilde{x}_1^2 & \tilde{x}_1 \tilde{x}_2\\
\tilde{x}_1 \tilde{x}_2 & \tilde{x}_2^2
\end{pmatrix}
  \bigg) \exp\big( - |\tilde{{\bm x}}|^2 \big).
  \end{align}

\section{Normalization and simulations}
\label{app:simul}

We performed numerical simulations of the paraxial 
wave equation \blue{Eq.(1)} (main text) by normalizing the spatial variables with respect to the wavelength $\lambda$:
\begin{equation}
\label{eq:parax_norm}
i \partial_{z'} \psi_{z'}(x') = - {\alpha'}  \partial_{x'}^2 \psi_{z'} + 
V'(z',x') \psi_{z'},
\end{equation}
where $x'=x/\lambda$, $z'=z/\lambda$, $V'=\lambda V=\pi (n_o^2-n^2(z',x'))/n_o$, and $\alpha'=\alpha/\lambda=1/(4\pi n_o) \simeq 0.053$ with a reference refractive index of $n_o=1.5$. 
Accordingly, the normalized initial correlation length is $\rho_o'=\rho_o/\lambda$, and the normalized variance of the random potential is $\sigma'^2=\EE[V'^2] =\lambda^2 \sigma^2$. Note that the relevant parameters are invariant with respect to the normalization, $X_c=X_c'=\sigma'^{2/3} \ell_c'/\alpha'^{1/3}$, $X_o=X_o'=\sigma'^{2/3} \rho_o'/\alpha'^{1/3}$, and $z_c/\lambda=z_c'=1/(2\sigma'^{2/3} \alpha'^{2/3})$.

The normalized paraxial Eq.(\ref{eq:parax_norm}) is solved using a pseudo-spectral split-step method, with a frequency cutoff of the spectral grid $k_c'=2\pi$ (i.e., $k_c=2\pi/\lambda=k_o$ refers to the light wavenumber in dimensional units), so that the spatial discretization is $dx'=1/2$ (i.e., $dx=\lambda/2$ in dimensional units). 
In all simulations, the size of the spatial window, $T_{x'}$, is chosen to be much larger than $\ell_c'$. Typically, we take $T_{x'}/\ell_c' \simeq 40 $.
Each realization of the random processes $V'(x',z')$ and $\psi_{z'=0}(x')$ are defined in the spectral domain using Gaussian correlation functions characterized by $\ell_c'$ and $\rho_o'$ respectively. The results presented in the main part of the text are the results of the numerical simulations, averaged over 1000 realizations in cases: 1) of an initial plane wave and 2) of a coherent speckled field, corresponding to situation {\tt (c)}. In the case of a partially coherent speckled initial field, situation {\tt (pc)}, we perform 300 realizations of the potential $V'(x',z')$. For each of those realizations, we perform an average over 400 realizations of the initial field $\psi_{z'=0}(x')$. The different realizations are performed in parallel using HPC resources from DNUM CCUB (Centre de Calcul de l'Universit\'e de Bourgogne).


\begin{thebibliography}{9}
\bibitem{comm1}
This results from the fact that, for $z$ large, the correlation function has the form $ \exp\big( - { {\gamma}_2 z}  y^2 /2 \big)$,
where $ {\gamma}_2=- \partial_x^2 \gamma(0) = \sigma^2 \ell_c^{-1} \tilde{\gamma}_2$,
$\tilde{\gamma}_2 =- \partial_{\tilde{x}}^2 \tilde{\gamma}(0) $, is always positive because it is proportional to the power spectral density of the stationary process $ \partial_x {V}(z,0)$ at frequency $0$.
\end{thebibliography}
\end{document}